\documentclass[a4paper,12pt]{amsart}
\usepackage{amssymb}

\title{Free Field Realization of \\ Vertex Operators for level two modules of $U_q\bigl(\widehat{\mathfrak{sl}}(2)\bigr)$}
\author{Yuji Hara}
\address{Institute of Physics, Graduate School of Arts and Sciences, University of Tokyo, Tokyo 153, Japan}
\email{ss77070@komaba.ecc.u-tokyo.ac.jp}
\begin{document}
\maketitle
\begin{abstract}
Free field realization of vertex operators for level two modules of $U_q\bigl(\widehat{\mathfrak{sl}}(2)\bigr)$ are shown through the free field realization of the modules given by Idzumi in Ref.\cite{Idzumi1,Idzumi2}. We constructed types I and II vertex operators when the spin of the associated evaluation module is 1/2 and type II's for the spin 1.
\end{abstract}
 
\section{Introduction}
Vertex operators for the quantum affine algebra $U_q\bigl(\widehat{\mathfrak{sl}}(2)\bigr)$ have played essential roles in the algebraic analysis of solvable lattice models since the pioneering works of \cite{CMP151,JMMN,JM}. In these works which analyze the XXZ model, type I vertex operators are identified with half infinite transfer matrices as their representation-theoretical counter part and type II vertex operators are interpreted as particle creation operators. To perform concrete computation such as a trace of composition of vertex operators, we need free field realization of modules and operators. In the said example of the XXZ model, the integral expressions of n-point correlation functions which are special cases of the traces are obtained through bosonization of level one module of $U_q\bigl(\widehat{\mathfrak{sl}}(2)\bigr)$.

Motivated by these results, Idzumi \cite{Idzumi1, Idzumi2} constructed level two modules and type I vertex operators accompanied by spin 1 evaluation modules for $U_q\bigl(\widehat{\mathfrak{sl}}(2)\bigr)$ in terms of bosons and fermions and then calculated correlation functions of a spin 1 analogue of the XXZ model.  
The purpose of this paper is to extend Idzumi's free field realization to other kinds of vertex operators i.e.\ type I and II vertex operators for the level two modules associated with the evalution module of spin 1/2 and the type II's for the spin 1. The results are given in Section \ref{sec:VOff} and their derivation is discussed in the first case in Section \ref{sec:derivation}. The results together with Ref.\cite{Idzumi1, Idzumi2} give the complete set of vertex operators for level two module of $U_q\bigl(\widehat{\mathfrak{sl}}(2)\bigr)$ and enable one to calculate form factors of the spin 1 analogue of the XXZ model.

Recently Jimbo and Shiraishi \cite{JS} showed a coset-type construction for the deformed Virasoro algebra with the vertex operators for $U_q\bigl(\widehat{\mathfrak{sl}}(2)\bigr)$. They constructed a primary operator for the deformed Virasoro algebra as coset type composition of vertex operators which may be denoted as $\Bigr(U_q\bigl(\widehat{\mathfrak{sl}}(2)\bigr)_k\oplus U_q\bigl(\widehat{\mathfrak{sl}}(2)\bigr)_1\Bigl)/ U_q\bigl(\widehat{\mathfrak{sl}}(2)\bigr)_{k+1}$. We hope that our results will be helpful for extending this work to the deformed supersymmetric Virasoro algebra through $\Bigr(U_q\bigl(\widehat{\mathfrak{sl}}(2)\bigr)_k\oplus U_q\bigl(\widehat{\mathfrak{sl}}(2)\bigr)_2\Bigl)/ U_q\bigl(\widehat{\mathfrak{sl}}(2)\bigr)_{k+2}$. 

\section{Free field realization of level two module}
\label{sec:convention}

\subsection{Convention}
In the following we will use $U$ to denote the quantum affine algebra $U_q\bigl(\widehat{\mathfrak{sl}}(2)\bigr)$.
Unless mentioned, we follow the notations of Ref.\cite{Idzumi1,Idzumi2}. As for the free field representation, we slightly modify the convention.

The quantum affine algebra $U$ is an associative algebra with unit 1 generated by $e_i,f_i\;(i=0,1),q^h\;(h\in P^*)$ with relations
\[q^0=1,\;q^h q^{h'}=q^{h+h'},\]
\[q^h e_i q^{-h}= q^{\langle h,\alpha_i\rangle}e_i,\;q^h f_i q^{-h}= q^{-\langle h,\alpha_i\rangle}f_i ,\]
\[\left[e_i,f_i\right]\footnotemark[1]=\delta_{ij}\frac{t_i-t_i^{-1}}{q-q^{-1}},\;(t_i=q^{h_i})\]
\[e_i^3e_j-[3]e_i^2e_je_i+[3]e_ie_je_i^2-e_je_i^3=0 ,\]
\[f_i^3f_j-[3]f_i^2f_jf_i+[3]f_if_jf_i^2-f_jf_i^3=0 ,\]
\footnotetext[1]{[A,B]=AB-BA}
\noindent where $P={\mathbb Z}\Lambda_0+{\mathbb Z}\Lambda_1+{\mathbb Z}\delta$ is the weight lattice of the affine Lie algebra $\widehat{\mathfrak{sl}}(2)$ and $P^*$ is the dual lattice to $P$ with the dual basis $\{h_0,h_1,d\}$ to $\{\Lambda_0,\Lambda_1,\delta\}$ with respect to the natural pairing $\langle\;\:,\;\rangle :P\times P^* \rightarrow {\mathbb Z}$. We also use current type generators introduced by Drinfeld \cite{Drin}
%\[a_k\;(k\in{\mathbb Z}\0),x^{\pm}_k\;(k\in{\mathbb Z}),\gamma{\pm1/2},K^\pm.\]\[\gamma:central\]
\[ [a_k.a_l]=\delta_{k+l,0}\frac{[2k]}{k}\frac{\gamma^k-\gamma^{-k}}{q-q^{-1}},\]
\[ Ka_kK^{-1}=a_k,\;\:Kx^\pm_kK^{-1}=q^{\pm 2}x^\pm_k,\]
\[ [a_k,x^\pm_l]=\pm\frac{[2k]}{k}\gamma^{\mp |k|/2}x^\pm_{k+l},\]
\[ x^\pm_{k+l}x^\pm_l-q^{\pm 2}x^\pm_lx^\pm_{k+l}=q^{\pm 2}x^\pm_kx^\pm_{l+1}-x^\pm_{l+1}x^\pm_k,\]
\[ [x^+_k,x^-_l]=\frac{\gamma^{\frac{k-l}{2}}\psi_{k+l}-\gamma^{\frac{l-k}{2}}\phi_{k+l}}{q-q^{-1}},\]
\noindent where $\psi_k,$ and $\varphi_k$ are defined as
\[\sum_{k\geqslant 0}\psi_k z^{-k}=K\exp \bigl\{(q-q^{-1})\sum_{k\geqslant 1}a_k z^{-k}\bigr\},\]
\[\sum_{k\geqslant 0}\phi_k z^{k}=K^{-1}\exp \bigl\{-(q-q^{-1})\sum_{k\geqslant 1}a_{-k} z^{k}\bigr\}.\]
\noindent The relation between two types of generators are
\[t_1=K,\;t_0=\gamma K^{-1},\;e_1=x^+_0,\;e_0t_1=x^-_1,\;f_1=x_0^-,\;t_1^{-1}f_1=x_0^{-1}.\]

The higest weight module and the evaluation module are described compactly in Ref.\cite{Idzumi1}.

Commutation and anti commutation relations of bosons and fermions are given by
\begin{gather*}
[a_m,a_n]=\delta_{m+n,0}\frac{[2m]^2}{m},\\
\{\phi_m,\phi_n\}\footnotemark[1]=\delta_{m+n,0}\eta_m,\\
\eta_m=q^{2m}+q^{-2m}.
\end{gather*}
\footnotetext[1]{\{A,B\}=AB-BA}
where $m,n\in\mathbb{Z}+1/2\mbox{ or }\in\mathbb{Z}$  for Neveu-Schwarz-sector or Ramond-sector respectively. Fock spaces and vacuum vectors are denoted as $\mathcal{F}^a,\;\mathcal{F}^{\phi^{NS}},\;\mathcal{F}^{\phi^{R}}\mbox{ and }|vac\rangle$,$|NS\rangle,|R\rangle$ for the boson and $NS$ and $R$ fermion respectively. Fermion currents are defined as
\begin{equation*}
\phi^{NS}(z)=\sum_{n\in\mathbb{Z}+\frac12}\phi^{NS}_nz^{-n},\;\phi^R(z)=\sum_{n\in\mathbb{Z}}\phi^R_nz^{-n}.
\end{equation*}

$Q=\mathbb{Z}\alpha$ is the root lattice of $\mathfrak{sl}_2$ and $F[Q]$ be the group algebra. We use $\partial$ as
\[
[\partial, \alpha]=2.
\]

\subsection{$V(2\Lambda_0),V(2\Lambda_1)$}
The highest weight module $V(2\Lambda_0)$ is identified with the Fock space
\begin{equation}
\mathcal{F}^{(0)}_+=\mathcal{F}^a\otimes\big{\{}(\mathcal{F}^{\phi^{NS}}_{even}
\otimes F[2Q])\oplus(\mathcal{F}^{\phi^{NS}}_{odd}\otimes e^{\alpha}F[2Q])\big{\}},
\end{equation}	
subscripts $even\mbox{ and }odd$ represent the number of fermions. The highest weight vector is $|vac\rangle\otimes |NS\rangle\otimes 1$.
$V(2\Lambda_1)$ is
\begin{equation}
\mathcal{F}^{(0)}_-=\mathcal{F}^a\otimes\big{\{}(\mathcal{F}^{\phi^{NS}}_{even}
\otimes e^\alpha F[2Q])\oplus(\mathcal{F}^{\phi^{NS}}_{odd}\otimes F[2Q])\big{\}}
\end{equation}	
with the highest weight vector being $|vac\rangle\otimes |NS\rangle\otimes e^\alpha$. Note that
\begin{gather*}
\mathcal{F}^{(0)}=\mathcal{F}^{(0)}_-\oplus\mathcal{F}^{(0)}_+.\nonumber\\
\mathcal{F}^{(0)}=\mathcal{F}^a\otimes\mathcal{F}^{\phi^{NS}}\otimes F[Q]\nonumber.
\end{gather*}
The operators are realized in the following manner.
\begin{gather*}
\gamma=q^2,\quad K=q^{\partial},\\
x^{\pm}(z)=\sum_{m\in\mathbb{Z}}x^{\pm}_mz^{-m}=E^\pm_<(z)E^\pm_>(z)\phi^{NS}(z)e^{\pm\alpha}z^{\frac12\pm\frac12\partial},
\end{gather*}
\begin{gather}
E^\pm_<(z)=\exp{(\pm\sum_{m>0}\frac{a_{-m}}{[2m]}q^{\mp m}z^{m})},\;\:E^\pm_>(z)=\exp{(\mp\sum_{m>0}\frac{a_m}{[2m]}q^{\mp m}z^{-m})},\nonumber\\ 
\intertext{and}
d=-\frac{\partial^2}{8}+\frac{(\lambda,\lambda)}{4}-\sum_{m=1}^\infty mN_m^a-\sum_{k>0}kN^{\phi^{NS}}_k,\label{eqn:d}\\
N^a_m=\frac{m}{[2m]^2}a_{-m}a_m,\quad N^{\phi^{NS}}_k=\frac{1}{\eta_m}\phi^{NS}_{-m}\phi^{NS}_m\quad (m>0),
\end{gather}
where the higest weight vector of the module should be substituted for $\lambda$ of \eqref{eqn:d}. 

\subsection{$V(\Lambda_0+\Lambda_1)$}
The module $V(\Lambda_0+\Lambda_1)$ is identified with
\begin{equation}
\mathcal{F}^{(1)}=\mathcal{F}^a\otimes\mathcal{F}^{\phi^{R}}\otimes 
e^{\frac{\alpha}{2}}F[Q],
\end{equation}
where 
\[
\phi^R_0 |R\rangle=|R\rangle.
\]
The highest weight vector is identified with $|vac\rangle\otimes |R\rangle\otimes e^{\frac{\alpha}{2}}$.

Operators are constructed in the same way as before except that subscripts for fermion sector are $R$ instead of $NS$.

\section{Free field realizations of vertex operators}
\label{sec:VOff}
Let $V,V'$ be level two modules and $V^{(k)}_z$ be a spin k/2 evaluation module of $U$. Vertex operators we will consider are $U$-linear maps of the following kinds \cite{FR,CMP155}
\begin{align}
\Phi_V^{V',k}(z)&:V\longrightarrow V'\otimes V_z^{(k)},\label{eqn:VOI}\\
\Psi_V^{k,V'}(z)&:V\longrightarrow V_z^{(k)}\otimes V'.\label{eqn:VOII}
\end{align}
Vertex operators of the form (\ref{eqn:VOI},\ref{eqn:VOII}) are called type I and II respectively. Components of vertex operators are defined as
%\newpage
\[
\Phi(z)_V^{V',k}=\sum_{n=0}^k\Phi_n(z)\otimes u_n,\;
\Psi(z)_V^{k,V'}=\sum_{n=0}^ku_n\otimes \Psi_n(z).
\]

\subsection{type I Vertex Operators for level 2 and spin 1/2}
We show free field realization of type I vertex operators of the following kinds\begin{align}
\Phi_{2\Lambda_i}^{\Lambda_0+\Lambda_1,1}(z)&:V(2\Lambda_i)\longrightarrow V(\Lambda_0+\Lambda_1)\otimes V_z^{(1)},\label{eqn:VOI1-1}\\
\Phi_{\Lambda_0+\Lambda_1}^{2\Lambda_i,1}(z)&:V(\Lambda_0+\Lambda_1)\longrightarrow V(2\Lambda_i)\otimes V_z^{(1)}\label{eqn:VOI1-2}
\end{align}
where $i=0\mbox{ or }1$. 

Under the free field realization of level 2 modules reviewed in Secton \ref{sec:convention}, the explicit forms of the components of the vertex operators in (\ref{eqn:VOI1-1}) are
\begin{align}
\Phi_1&(z)=B_{I,<}(z)B_{I,>}(z)
        \Omega_{NS}^R(z)e^{\alpha/2}(-q^4z)^{\partial/4},\label{eqn:VOI1-ff1}\\
\Phi_0&(z)=\oint\frac{dw}{2\pi i}B_{I,<}(z)E^-_<(w)B_{I,>}(z)E^-_>(w)
    \Omega_{NS}^R(z)\phi^{NS}(w)\label{eqn:VOI1-ff2}\\
    &\times e^{-\alpha/2}(-q^4z)^{\partial/4}w^{-\partial/2}
    (-q^4zw^3)^{-\frac 12}\frac{\bigr(\displaystyle{\frac{w}{q^3z}};q^4\bigl)_{\infty}}{\bigr(\displaystyle{\frac{w}{qz}};q^4\bigl)_{\infty}}\Bigl\{\frac{w}{1-q^{-3}w/z}+\:\frac{q^5z}{1-q^5z/w}\Bigr\}\nonumber,
\end{align}
\begin{gather}
B_{I,<}(z)= \exp\bigl(\sum_{n=1}^{\infty}\frac{[n]a_{-n}}{[2n]^2}(q^5z)^n\bigr),\\
B_{I,>}(z)= \exp\bigl(-\sum_{n=1}^{\infty}\frac{[n]a_{n}}{[2n]^2}(q^3z)^{-n}\bigr).
\end{gather}
The integrand of $\Phi_0(z)$ has poles only at $w=q^5z,q^3z$ except for $w=0, \infty$ and the contour of integration encloses $w=0,q^5z$, details are discussed in Sec.\ref{sec:derivation}.
For those of (\ref{eqn:VOI1-2}) we just replace $\Omega_{NS}^R(z)$ with $\Omega_{R}^{NS}(z)$ in (\ref{eqn:VOI1-ff1},\ref{eqn:VOI1-ff2}).

The fermionic part $\Omega(z)$'s are maps between different fermion sectors and satisfy 
\begin{equation}
\phi^{NS}(w)\Omega(z)^{NS}_R=\Bigr(\frac{-q^4z}{w}\Bigl)^{1/2}\frac{\bigr(\displaystyle{\frac{w}{q^3z}};q^4\bigl)_\infty\bigr(\displaystyle{\frac{q^7z}{w}};q^4\bigl)_\infty}{\bigr(\displaystyle{\frac{w}{qz}};q^4\bigl)_\infty\bigr(\displaystyle{\frac{q^5z}{w}};q^4\bigl)_\infty}\Omega(z)^{NS}_R\phi^R(w)
\label{eqn:FE-comm}
\end{equation}
and exactly the same equation except subscripts for fermion sectors are exchanged.
This kind of mapping for fermions first appeared in high-energy phisics theory as \lq\lq fermion emission vertex operator"\cite{JMP35,FE}. Their free field realizations are
\begin{gather}
\Omega_{NS}^R(z)=\langle NS|\ e^Y|R\rangle,\label{eqn:FE-1}\\[10pt]
\begin{split}
Y&=-\sum_{m>n\geq 0}X_{m,n}\varphi_{-m}^{R}\varphi_{-n}^{R}z^{m+n}
  -\sum_{k>l\geq 0}X_{k+1/2,l+1/2}\varphi_{k+1/2}^{NS}\varphi_{l+1/2}^{NS}z^{-k-l-1}\\
  &+\sum_{
\begin{subarray}{c}
m\geq 0\\ k\geq 0
\end{subarray}}
X_{m,-k-1/2}\varphi_{-m}^{R}\varphi_{k+1/2}^{NS}z^{m-k-1/2}
\end{split}\label{eqn:FE-2}
\end{gather}
\begin{gather}
\Omega_{R}^{NS}(z)=\langle R|\ e^{Y'}|NS\rangle,\label{eqn:FE-3}\\[10pt]
\begin{split}
Y'&=\sum_{k>l\geq 0}X_{k+1/2,l+1/2}\varphi_{-k-1/2}^{NS}\varphi_{-l-1/2}^{NS}z^{k+l+1}   +\sum_{m>n\geq 0}X_{m,n}\varphi_{m}^{R}\varphi_{n}^{R}z^{-m-n}\\
  &-\sum_{
\begin{subarray}{c}
k\geq 0\\ m\geq 0
\end{subarray}}
X_{-k-1/2,m}\varphi_{-k-1/2}^{NS}\varphi_{m}^{R}z^{k-m+1/2}
\end{split}
\end{gather}
\begin{gather}
\varphi_0^R=\phi_0^R,\quad
\varphi_{-m}^{R}=\phi_{-m}^{R}\frac{\gamma_m q^{5m}}{\eta_m},\quad
\varphi_{m}^{R}=\phi_{m}^{R}\frac{\gamma_m q^{-3m}}{\eta_m}\quad (m>0),\label{eqn:varphi-1}\\
\varphi_{k+1/2}^{NS}=\phi_{k+1/2}^{NS}\frac{\gamma_k q^{-3k-2}}{\eta_{k+1/2}}(-(-1)^{1/2}),\quad\varphi_{-k-1/2}^{NS}=\phi_{-k-1/2}^{NS}\frac{\gamma_k q^{5k+2}}{\eta_{k+1/2}}(-1)^{1/2}\quad (k>0),\label{eqn:varphi-2}\\
X_{k,l}=\frac{q^{4k}-q^{4l}}{1-q^{4(k+l)}},\nonumber\\
\gamma_n=\frac{(q^2;q^4)_n}{(q^4;q^4)_n},\quad \frac{(q^2z;q^4)_\infty}{(z;q^4)_\infty}=\sum^\infty_{n=0}\gamma_n z^n\label{eqn:q-binomial}.
\end{gather}
(\ref{eqn:FE-1},\ref{eqn:FE-3}) are to mean that a matrix element is given by
\begin{gather*}
{}_R\langle\mbox{out}|\Omega_{NS}^R(z)|\mbox{in}\rangle_{NS}={}_R\langle\mbox{out}|\otimes\langle NS|\ e^{Y}|R\rangle\otimes|\mbox{in}\rangle_{NS},\\
\mbox{for   }
|\mbox{out}\rangle_R\in\mathcal{F}^{\phi^{R}},\;|\mbox{in}\rangle_{NS}\in\mathcal{F}^{\phi^{NS}}.
\end{gather*}

We define the normalized vertex operators $\tilde{\Phi}(z)$'s as follows
\begin{gather*}
\langle\Lambda_0+\Lambda_1|\tilde{\Phi}_1(z)|2\Lambda_0\rangle=1,\quad
\langle 2\Lambda_1|\tilde{\Phi}_1(z)|\Lambda_0+\Lambda_1\rangle=1,\\
\langle\Lambda_0+\Lambda_1|\tilde{\Phi}_0(z)|2\Lambda_1\rangle=1,\quad
\langle 2\Lambda_0|\tilde{\Phi}_0(z)|\Lambda_0+\Lambda_1\rangle=1,
\end{gather*}
and these are given by
\begin{gather}
\tilde{\Phi}_{2\Lambda_0}^{\Lambda_0+\Lambda_1,1}(z)=\Phi(z),\\
\tilde{\Phi}^{2\Lambda_1,1}_{\Lambda_0+\Lambda_1}(z)=(-q^4z)^{-1/4}\Phi(z),\\
\tilde{\Phi}^{2\Lambda_0,1}_{\Lambda_0+\Lambda_1}(z)=(-q^4z)^{1/4}\Phi(z),\\
\tilde{\Phi}_{2\Lambda_1}^{\Lambda_0+\Lambda_1,1}(z)=(-q^6z)^{-1/2}\Phi(z).
\end{gather}

\subsection{type II Vertex Operators for level 2 and spin 1/2}
We consider type II vertex operators of the following kind
\begin{align}
\Psi_{2\Lambda_i}^{1,\Lambda_0+\Lambda_1}(z)&:V(2\Lambda_i)\longrightarrow V_z^{(1)}\otimes V(\Lambda_0+\Lambda_1),\\
\Psi_{\Lambda_0+\Lambda_1}^{1,2\Lambda_i}(z)&:V(\Lambda_0+\Lambda_1)\longrightarrow V_z^{(1)}\otimes V(2\Lambda_i).
\end{align}
Explicit forms of the components are as follows.
\begin{gather}
\begin{align}
\Psi_0(z)&=B_{II,<}(z)B_{II,>}(z)
         \Omega(q^{-2}z)e^{-\alpha/2}(-q^2z)^{-\partial/4},\\
\Psi_1(z)&=\oint\frac{dw}{2\pi i}B_{II,<}(z)E^+_<(w)B_{II,>}(z)E^+_>(w)
    \Omega(q^{-2}z)\phi(w)\\
    &\times e^{\alpha/2}(-q^2z)^{-\partial/4}w^{\partial/2}
    (-q^2zw^3)^{-\frac 12}\frac{\bigr(\displaystyle{\frac {w}{qz}};q^4\bigl)_{\infty}}{\bigr(\displaystyle{\frac{qw}{z}};q^4\bigl)_{\infty}}\Bigl\{\frac{w}{1-q^{-3}w/z}+\frac{q^3z}{1-qz/w}\Bigr\}\nonumber,
\end{align}\\
B_{II,<}(z)= \exp\bigl(-\sum_{n=1}^{\infty}\frac{[n]a_{-n}}{[2n]^2}(qz)^n\bigr),\\
B_{II,>}(z)= \exp\bigl(\sum_{n=1}^{\infty}\frac{[n]a_{n}}{[2n]^2}(q^3z)^{-n}\bigr).
\end{gather}
The integrand of $\Psi_1(z)$ has poles only at $w=q^3z,qz$ except for $w=0, \infty$ and the contour of integration encloses $w=0,qz$.
Subscripts for fermion sectors are abbreviated.

Normalized vertex operators are defined by the conditions 
\begin{gather*}
\langle\Lambda_0+\Lambda_1|\tilde{\Psi}_1(z)|2\Lambda_0\rangle=1,\quad
\langle 2\Lambda_1|\tilde{\Psi}_1(z)|\Lambda_0+\Lambda_1\rangle=1,\\
\langle\Lambda_0+\Lambda_1|\tilde{\Psi}_0(z)|2\Lambda_1\rangle=1,\quad
\langle 2\Lambda_0|\tilde{\Psi}_0(z)|\Lambda_0+\Lambda_1\rangle=1,
\end{gather*}
and these are given by
\begin{gather}
\tilde{\Psi}_{2\Lambda_0}^{1,\Lambda_0+\Lambda_1}(z)=(-q)^{-1}\Psi(z),\\
\tilde{\Psi}^{1,2\Lambda_1}_{\Lambda_0+\Lambda_1}(z)=-(-q^6z)^{-1/4}\Psi(z),\\
\tilde{\Psi}^{1,2\Lambda_0}_{\Lambda_0+\Lambda_1}(z)=(-q^2z)^{1/4}\Psi(z),\\
\tilde{\Psi}_{2\Lambda_1}^{1,\Lambda_0+\Lambda_1}(z)=(-q^2z)^{1/2}\Psi(z).
\end{gather}

\subsection{type II Vertex Operators for level 2 and spin 1}
When the spin of the evaluation module is 1, the type II vertex operators do not contain any fermion emission vertex operators. 
\begin{align}
\Psi_{2\Lambda_i}^{2,2\Lambda_i}(z)&:V(2\Lambda_i)\longrightarrow V_z^{(2)}\otimes V(2\Lambda_i),\\
\Psi_{\Lambda_0+\Lambda_1}^{2,\Lambda_0+\Lambda_1}(z)&:V(\Lambda_0+\Lambda_1)\longrightarrow V_z^{(2)}\otimes V(\Lambda_0+\Lambda_1).
\end{align}
Explicit form of the components are as follows.
\begin{align}
\Psi_0(z)&=F_{II,<}(z)F_{II,>}(z)e^{-\alpha}(-q^2z)^{-\partial/2+1},\\
\Psi_1(z)&=\oint\frac{dw}{2\pi i}F_{II,<}(z)E^+_<(w)F_{II,>}(z)E^+_>(w)\phi(w)\Bigr(\frac{w}{-q^2z}\Bigl)^{\partial/2}\\
&\times w^{-1/2}\Bigr\{\frac{1}{1-\displaystyle{\frac{w}{q^4z}}}+\frac{q^4z}{w\bigr(1-\displaystyle{\frac{z}{w}}\bigl)}\Bigl\}\nonumber,
\end{align}
The integration contour encircles poles $w=0,z$ but the pole $w=q^4z$ lies outside of it.
\begin{gather}
\begin{align}
\Psi_2(z)=\oint\frac{dw_2}{2\pi i}&\oint\frac{dw_1}{2\pi i}F_{II,<}(z)E^+_<(w_1)E^+_<(w_2)F_{II,>}(z)E^+_>(w_1)E^+_>(w_2)\\ 
&\times e^{\alpha}\Bigr(\frac{w_1w_2}{-q^2z}\Bigl)^{\partial/2}(w_1w_2)^{-1/2}\Bigr\{\frac{1}{1-\displaystyle{\frac{w_1}{q^4z}}}+\frac{q^4z}{w_1\bigr(1-\displaystyle{\frac{z}{w_1}}\bigl)}\Bigl\}\nonumber\\ &\times\Bigr\{[2]^{-1}:\phi(w_1)\phi(w_2):\Biggr(\frac{w_1-q^{-2}w_2}{-q^2z\bigr(1-\displaystyle{\frac{w_2}{q^4w_1}}\bigl)}+\frac{1-\displaystyle{\frac{w_1}{q^2w_2}}}{1-\displaystyle{\frac{z}{w_2}}}\Biggl)\nonumber\\ 
&\quad\quad+\frac{(w_1w_2)^{1/2}\bigr(1-\displaystyle{\frac{w_2}{w_1}}\bigl)}{-q^2z\bigr(1-\displaystyle{\frac{q^2w_2}{w_1}}\bigl)\bigr(1-\displaystyle{\frac{w_2}{q^4z}}\bigl)}-\frac{\displaystyle{\bigr(\frac{w_1}{w_2}\bigl)^{1/2}}\bigr(1-\displaystyle{\frac{w_1}{w_2}}\bigl)}{\bigr(1-\displaystyle{\frac{q^2w_1}{w_2}})\bigr(1-\displaystyle{\frac{z}{w_2}}\bigl)}\Bigl\},\nonumber
\end{align}
\end{gather}
We have to prepare two contours because of the fermionic part and one is for the term including $:\phi(w_1)\phi(w_2):$ and the other is for the rest. The former satisfies $|\frac{w_2}{q^4w_1}|<1,|w_2|>|z|$ and the same condition satisfied by the contour for $Psi_1$ with substitution $w=w_1$. The latter satisfies $|q^2w_2|<|w_1|<|q^{-2}w_2|$ and the same conditions as $Psi_1$ with $w=w_1,w_2$.
\begin{gather}
F_{II,<}(z)=\exp\Bigr(-\sum_{m>0}\frac{a_{-m}}{[2m]}(qz)^m\Bigl),\\
F_{II,>}(z)=\exp\Bigr(\sum_{m>0}\frac{a_{m}}{[2m]}(q^3z)^{-m}\Bigl).
\end{gather}
Under the normailzation
\begin{gather*}
\langle 2\Lambda_0|\tilde{\Psi}_0(z)|2\Lambda_1\rangle=1,\quad
\langle 2\Lambda_1|\tilde{\Psi}_2(z)|2\Lambda_0\rangle=1,\\
\langle\Lambda_0+\Lambda_1|\tilde{\Psi}_1(z)|\Lambda_0+\Lambda_1\rangle=1,
\end{gather*}
$\tilde{\Psi}(z)$' are given by
\begin{gather}
\tilde{\Psi}_{2\Lambda_1}^{2,2\Lambda_0}(z)=\Psi(z),\\
\tilde{\Psi}_{\Lambda_0+\Lambda_1}^{2,\Lambda_0+\Lambda_1}(z)=-(-q^2z)^{-1/2}\Psi(z),\\
\tilde{\Psi}_{2\Lambda_0}^{2,2\Lambda_1}(z)=(-q^4z)^{-1}\Psi(z).
\end{gather}

\section{Derivation}
\label{sec:derivation}
Taking $\Phi_{2\Lambda_i}^{\Lambda_0+\Lambda_1,1}(z)$ as an example, we discuss the derivation of the results in the previous section. Other cases can be treated in almost the same way.

\subsection{General structure of $\Phi_0(z)\mbox{ and }\Phi_1(z)$}
Calculating
\[
\Delta(x)\Phi(z)=\Phi(z)x
\]
\newpage
for $x=\mbox{Chevalley generators of $U$ and }a_n$, we get
\begin{gather}
0=[\Phi_1(z),x_0^+],\nonumber\\
K\Phi_1(z)=[\Phi_0(z),x_0^+],\nonumber\\
0=x_0^-\Phi_0(z)-q\Phi_0(z)x_0^-,\nonumber\\
\Phi_0(z)=\Phi_1(z)x_0^--qx_0^-\Phi_1(z),\label{eqn:Chv-2}\\
0=\Phi_0(z)x_1^--qx_1^-\Phi_0(z),\nonumber\\
q^{3}z\Phi_0(z)=\Phi_1(z)x_{1}^--q^{-1}x_{1}^-\Phi_1(z),\label{eqn:Chv-3}\\
(qzK)^{-1}\Phi_1(z)=[\Phi_0(z),x_{-1}^+],\nonumber\\
0=[\Phi_1(z),x_{-1}^+],\nonumber\\
K\Phi_1(z)K^{-1}=q\Phi_1(z),\label{eqn:Chv-1}\\
K\Phi_0(z)K^{-1}=q^{-1}\Phi_0(z),\nonumber\\
%\end{gather}
%\begin{gather}
[a_m,\Phi_1(z)]=(q^5z)^m\frac{[m]}{m}\Phi_1(z)\label{eqn:boson-1},\\
[a_{-m},\Phi_1(z)]=(q^3z)^{-m}\frac{[m]}{m}\Phi_1(z)\label{eqn:boson-2}.
\end{gather}
From (\ref{eqn:Chv-1},\ref{eqn:boson-1},\ref{eqn:boson-2}), we can speculate the form of $\Phi_1(z)$ as
\[
\Phi_1(z)=B_{I,<}(z)B_{I,>}(z)
         \Omega^R_{NS}(z)e^{\alpha/2}y^{\partial}.
\]
To determine $y$ and the fermionic part $\Omega^R_{NS}(z)$, we impose the following conditions on $\Phi_1(z)$
\begin{gather*}
\Phi_1(z)x_0^--qx_0^-\Phi_1(z)=(q^3z)^{-1}(\Phi_1(z)x_{1}^--q^{-1}x_{1}^-\Phi_1(z)),\\
0=[\Phi_1(z),x^+(w)],
\end{gather*}
which can be easily seen from (\ref{eqn:Chv-2},\ref{eqn:Chv-3}) and the proposition of \ Section 4.4 of Ref.\cite{ChPr}. Then we have (\ref{eqn:VOI1-ff1},\ref{eqn:FE-comm})
\begin{gather*}
\Phi_1(z)=B_{I,<}(z)B_{I,>}(z)
    \Omega_{NS}^R(z)e^{\alpha/2}(-q^4z)^{\partial/4},\\
\phi^R(w)\Omega_{NS}^R(z)=\Bigr(\frac{-q^4z}{w}\Bigl)^{1/2}\frac{\bigr(\displaystyle{\frac{w}{q^3z}};q^4\bigl)_\infty\bigr(\displaystyle{\frac{q^7z}{w}};q^4\bigl)_\infty}{\bigr(\displaystyle{\frac{w}{qz}};q^4\bigl)_\infty\bigr(\displaystyle{\frac{q^5z}{w}};q^4\bigl)_\infty}\Omega_{NS}^R(z)\phi^{NS}(w).
\end{gather*}
$\Phi_1(z)$ can be calculated through (\ref{eqn:Chv-2})
\begin{align*}
\Phi_0(z)
    &=\oint\frac{dw}{2\pi i}\frac 1w\{\Phi_1(z)x^-(w)-qx^-(w)\Phi_1(z)\}\\
    &=\oint\frac{dw}{2\pi i}B_{I,<}(z)E^-_<(w)B_{I,>}(z)E^-_>(w)
    \Omega_{NS}^R(z)\phi^{NS}(w)\label{eqn:VOI1-ff2}\\
    &\times e^{-\alpha/2}(-q^4z)^{\partial/4}w^{-\partial/2}
    (-q^4zw^3)^{-\frac 12}\frac{\bigr(\displaystyle{\frac{w}{q^3z}};q^4\bigl)_{\infty}}{\bigr(\displaystyle{\frac{w}{qz}};q^4\bigl)_{\infty}}\Bigl\{\frac{w}{1-q^{-3}w/z}+\:\frac{q^5z}{1-q^5z/w}\Bigr\}\nonumber,
\end{align*}
To determine the contour of integration we have to find the poles of $\Omega_{NS}^R(z)\phi^{NS}(w)$ and this can be seen from
\[
\langle R|\Omega_{NS}^R(z)\phi^{NS}(w)|NS\rangle=\frac{{\bigr(\displaystyle{\frac{w}{qz}};q^4\bigl)_{\infty}}}{{\bigr(\displaystyle{\frac{w}{q^3z}};q^4\bigl)_{\infty}}},\quad \langle NS|\Omega^{NS}_R(z)\phi^{R}(w)|R\rangle=\bigl(\frac{w}{-q^4z}\bigr)^{1/2}\frac{{\bigr(\displaystyle{\frac{w}{qz}};q^4\bigl)_{\infty}}}{{\bigr(\displaystyle{\frac{w}{q^3z}};q^4\bigl)_{\infty}}}.
\]
Hence as a composite $\Omega_{NS}^R(z)\phi^{NS}(w)\frac{\bigr(\displaystyle{\frac{w}{q^3z}};q^4\bigl)_{\infty}}{\bigr(\displaystyle{\frac{w}{qz}};q^4\bigl)_{\infty}}$  in the integrand has no poles and the contour is the one encloses $w=0,q^5z$.

\subsection{Fermion emission vertex operator}
In Ref.\cite{JMP35}, Eqn.\eqref{eqn:FE-1} appears in the study of the Ising model and its free field realization is given without any details. Thus we give the exposition of its derivation\footnote[1]{We are indebted to M.Jimbo for explaining the details of Ref.\cite{JMP35}.}. The main point of derivating free field realization of the fermion emission vertex operator $\Omega^R_{NS}(z)$  (\ref{eqn:FE-1},\ref{eqn:FE-2}) is  to expand $\Omega^R_{NS}(z)$ as
\begin{gather*}
\Omega^R_{NS}(z)=\sum_{K,L}a_{K,L}\;\phi^R_{k_1}\phi^R_{k_2}\cdot\cdot\cdot|R\rangle\langle NS|\phi^{NS}_{l_1}\phi^{NS}_{l_2}\cdot\cdot\cdot ,\\
K=\{k_i\},L=\{l_i\},
\end{gather*}
and to calculate the coefficients $a_{K,L}$. After normalizing $\phi_n$ suitably to $\varphi_n$ (\ref{eqn:varphi-1},\ref{eqn:varphi-2}), we see \lq\lq $a_{K,L}/(\mbox{normalization factor})"$ are identified with Pfaffians of $X_{k,l}$. With the aid of a relation satisfied by Pfaffian
\begin{gather*}
\omega^{\wedge n}=n!\mbox{Pf}(b_{ij})x_1\wedge x_2\cdot\cdot\cdot\wedge x_{2n},\\
\intertext{where $x_k\:(1\leq k\leq 2n)$ is a Grassmann variable and}
\omega=\sum_{1\leq i<j\leq 2n}b_{ij}x_i\wedge x_j,
\end{gather*}
we get (\ref{eqn:FE-1},\ref{eqn:FE-2}).

Wick's theorem can be generalized to the present situation and we only need to calculate one- and two-point correlation functions for $a_{K,L}$. To calculate these, we rewrite \eqref{eqn:FE-comm} and introduce auxiliary operators
\begin{gather}
\tilde{\phi}^{NS}(w)\Omega^{NS}_{R}(q^{-4})=\Omega^{NS}_{R}(q^{-4})\tilde{\phi}^{R}(w),\label{eqn:FE1}\\
\tilde{\phi}^{NS}(w)=(-1)^{-1/2}w^{1/2}\frac{(qw^{-1};q^4)_\infty}{(q^3w^{-1};q^4)_\infty}\phi^{NS}(w),\\
\tilde{\phi}^{R}(w)=\frac{(qw;q^4)_\infty}{(q^3w;q^4)_\infty}\phi^{R}(w)=f_+(w)\phi^R(w)\label{eqn:aux-1},
\end{gather}
we set $\Omega(z=q^4)$ for simplicity. They are defined to satisfy 
\[
\langle {NS}|\tilde{\phi}^{NS}_n=0\;(n<0),\quad\tilde{\phi}^{R}_n|R\rangle=0\;(n>0),\quad\tilde{\phi}^{R}_0|R\rangle=|R\rangle,
\]
and this enables us to see that 
\[
\langle NS|\Omega_R^{NS}(q^{-4})\tilde{\phi}^R(z)\tilde{\phi}^R(w)|NS\rangle=
\langle NS|\tilde{\phi}^{NS}(z)\Omega_R^{NS}(q^{-4})\tilde{\phi}^R(w)|NS\rangle
\]
contains only negative (positive) powers of $z\;(w)$. On the other hand the expectation value of
\begin{gather*}
\{\tilde{\phi}^R(z),\tilde{\phi}^R(w)\}=f_+(z)f_+(w)\Bigl(\delta\bigl(\frac{q^2w}{z}\bigr)+\delta\bigl(\frac{w}{q^2z}\bigr)\Bigr),\\
\delta(z)=\sum_{n\in\mathbb{Z}}z^n,
\end{gather*}
with respect to $\langle NS|\Omega_R^{NS}(q^{-4})\mbox{ and }|R\rangle$ is
\begin{gather*}
\langle NS|\Omega_R^{NS}(q^{-4})\tilde{\phi}^R(z)\tilde{\phi}^R(w)|NS\rangle+
\langle NS|\Omega_R^{NS}(q^{-4})\tilde{\phi}^R(w)\tilde{\phi}^R(z)|NS\rangle\\
=f_+(z)f_+(w)\Bigl(\delta\bigl(\frac{q^2w}{z}\bigr)+\delta\bigl(\frac{w}{q^2z}\bigr)\Bigr)
\end{gather*}
where we normalize $\langle NS|\Omega_R^{NS}(q^{-4})|R\rangle=1$. And we get
\[
\langle NS|\Omega_R^{NS}(q^{-4})\tilde{\phi}^R(z)\tilde{\phi}^R(w)|R\rangle=\frac{1-qw}{1-q^2w/z}+\frac{1-q^{-1}w}{1-q^{-2}w/z}-1
\]
 Expanding the last line of the following equation as in \ref{sec:2pt}
\begin{align*}
\langle NS|\Omega_R^{NS}&(q^{-4})\phi^R(z)\phi^R(w)|R\rangle
  =\sum_{n,m\in\mathbb{Z}}\langle NS|\Omega_R^{NS}(q^{-4})\phi^R_n\phi^R_m|R\rangle z^{-n}w^{-m}\\
  &=\frac{1}{f_+(z)f_+(w)}\Bigl\{\frac{1-qw}{1-q^2w/z}+\frac{1-q^{-1}w}{1-q^{-2}w/z}-1\Bigr\},
\end{align*}
we have
\begin{equation}
\langle NS|\Omega^{NS}_R(q^{-4})\phi^R_{-n}\phi^R_{-m}|R\rangle=X_{m,n}\gamma_n\gamma_m q^{n+m}\;(n,m\geq 0).\label{eqn:2pt}
\end{equation}
Similar calculation yields
\begin{gather}
\langle NS|\phi^{NS}_{k+1/2}\Omega^{NS}_R(q^{-4})\phi^R_{-n}|R\rangle=-(-1)^{1/2}X_{-k-1/2,n}\gamma_n\gamma_k q^{n+k}\;(n,k\geq 0),\\[5pt]
\langle NS|\phi^{NS}_{k+1/2}\phi^{NS}_{l+1/2}\Omega^{NS}_R(q^{-4})|R\rangle=-X_{l+1/2,k+1/2}\gamma_l\gamma_k q^{l+k}\;(k,l\geq 0).\label{eqn:FE2}
\end{gather}

$z$-dependence of $\Omega^R_{NS}(z)$ is recovered with the equation
\begin{gather}
\zeta^{d^{R}}\Omega^R_{NS}(z)\zeta^{-d^{NS}}=\Omega^R_{NS}(\zeta^{-1}z),\label{eqn:z-dep}\\
\zeta^{-d^i}\phi^i(z)\zeta^{d^i}=\phi^i(\zeta z),\nonumber\\
\langle i|d^i=d^i|i\rangle=0,\nonumber
\end{gather}
where $d^i$'s are the fermionic part of $d$ of \eqref{eqn:d}
\begin{gather*}
d^i=-\sum_{k>0}kN^{\phi^{i}}_k,\;(i=NS\mbox{ or }R)\\
\intertext{and satisfy}
[d^i,\phi^i_n]=n\phi_n.
\end{gather*}
To derive \eqref{eqn:z-dep}, we multiply \eqref{eqn:FE-comm} by $\zeta^{d^{R}},\zeta^{-d^{NS}}$ from left and right respectively.

\bigskip
\noindent {\Large {\em Acknowledgement}}\\
The author thanks M.Jimbo, H.Konno, S.Odake and J.Shiraishi for helpful discussions. He also thanks A.Kuniba for warm encourragement.

\appendix

\section{boson}
Followings are useful formulae for normal ordering bosons. We set $(z)_\infty=(z;q^4)_\infty$ for brevity.
\begin{gather*}
B_{I,>}(z)E^-_<(w)=\frac{(qw/z)_\infty}{(q^{-1}w/z)_\infty}E^-_<(w)B_{I,>}(z),\\ E^-_>(w)B_{I,<}(z)=\frac{(q^9z/w)_\infty}{(q^7z/w)_\infty}B_{I,<}(z)E^-_>(w),\\
B_{I,>}(z)E^+_<(w)=\frac{(q^{-3}w/z)_\infty}{(q^{-1}w/z)_\infty}E^+_<(w)B_{I,>}(z),\\ E^+_>(w)B_{I,<}(z)=\frac{(q^5z/w)_\infty}{(q^7z/w)_\infty}B_{I,<}(z)E^+_>(w),\\
%\end{gather*}
%\begin{gather*}
B_{II,>}(z)E^+_<(w)=\frac{(q^{-1}w/z)_\infty}{(q^{-3}w/z)_\infty}E^+_<(w)B_{II,>}(z),\\ E^+_>(w)B_{II,<}(z)=\frac{(q^3z/w)_\infty}{(qz/w)_\infty}B_{II,<}(z)E^+_>(w), 
\end{gather*}
\begin{gather*}
B_{II,>}(z)E^-_<(w)=\frac{(q^{-1}w/z)_\infty}{(qw/z)_\infty}E^-_<(w)B_{II,>}(z),\\ E^-_>(w)B_{II,<}(z)=\frac{(q^3z/w)_\infty}{(q^5z/w)_\infty}B_{II,<}(z)E^-_>(w), \\
%\end{gather*}
%\begin{gather*}
F_{II,>}(z)E^-_<(w)=(1-\frac{w}{q^2z})E^-_<(w)F_{II,>}(z),\\ E^-_>(w)F_{II,<}(z)= (1-\frac{q^2z}{w})F_{II,<}(z)E^-_>(w),\\
F_{II,>}(z)E^+_<(w)=\frac{1}{1-q^{-4}w/z}E^+_<(w)F_{II,>}(z),\\ E^+_>(w)F_{II,<}(z)=\frac{1}{1-z/w}F_{II,<}(z)E^-_>(w),\\
E^-_>(w_1)E^+_<(w_2)=\frac{1}{1-w_2/w_1}E^+_<(w_2)E^-_>(w_1),\\ E^+_>(w_2)E^-_<(w_1)=\frac{1}{1-w_1/w_2}E^-_<(w_1)E^+_>(w_2), 
\end{gather*}

\section{fermion}
For $\Omega^{NS}_{R}(z)$, we show the equations corresponding to the ones from \eqref{eqn:FE1} to \eqref{eqn:FE2} 
\begin{gather}
\tilde{\phi}^{R'}(w)\Omega^R_{NS}(q^{-4})=\Omega^R_{NS}(q^{-4})\tilde{\phi}^{NS'}(w),\\
\tilde{\phi}^{R'}(w)=\frac{(q/w;q^4)_\infty}{(q^3/w;q^4)_\infty}\phi^R(w),\\
\tilde{\phi}^{NS'}(w)=(-1)^{1/2}w^{-1/2}\frac{(qw;q^4)_\infty}{(q^3w;q^4)_\infty}\phi^{NS}(w),
\end{gather}
\[
\langle R|\tilde{\phi}^{R'}_n=0\;(n<0),\quad\langle R|\tilde{\phi}^{R'}_0=\langle R|,\quad\tilde{\phi}^{NS'}_n|NS\rangle=0\;(n>0),
\]
\[
\langle R|\tilde{\phi}^{R'}(z)\tilde{\phi}^{R'}(w)\Omega^R_{NS}(q^{-4})|NS\rangle=\frac{1-q/z}{1-q^2w/z}+\frac{1-q^{-1}/z}{1-q^{-2}w/z}-1
\]
\begin{gather*}
\langle R|\phi^R_n\phi^R_m\Omega^R_{NS}(q^{-4})|NS\rangle=X_{n,m}\gamma_n\gamma_m q^{n+m}\;(n,m\geq 0),\\[5pt]
\langle R|\phi^R_n\Omega^R_{NS}(q^{-4})\phi^{NS}_{-k-1/2}|NS\rangle=(-1)^{1/2}X_{-k-1/2,n}\gamma_n\gamma_k q^{n+k}\;(n,k\geq 0),\\[5pt]
\langle R|\Omega^R_{NS}(q^{-4})\phi^{NS}_{-k-1/2}\phi^{NS}_{-l-1/2}|NS\rangle=X_{l+1/2,k+1/2}\gamma_l\gamma_k q^{l+k}\;(k,l\geq 0)
\end{gather*}

\section{Calculation of Eqn.\eqref{eqn:2pt}}
\label{sec:2pt}
We show details of calculation of \eqref{eqn:2pt}. From \eqref{eqn:q-binomial}
\begin{align*}
&\langle NS|\Omega_R^{NS}(q^{-4})\phi^R(z)\phi^R(w)|R\rangle\\
&=\frac{1}{f_+(z)f_+(w)}\Bigl\{\frac{1-qw}{1-q^2w/z}+\frac{1-q^{-1}w}{1-q^{-2}w/z}-1\Bigr\}\\
&=\sum_{k\geqslant 0,l\geqslant 0}\gamma_k (qz)^k\gamma_l (qw)^l\Bigl\{\sum_{a\geqslant 0}\Bigl(  (1-qw)\bigl(\frac{q^2w}{z}\bigr)^a+(1-w/q)\bigl(\frac{w}{q^2z}\bigr)^a  \Bigr)-1\Bigr\}\\
&=\sum_{0\leqslant a\leqslant m}\gamma_{n+a}\gamma_{m-a}\eta_a q^{n+m}z^n w^m-\sum_{0\leqslant a\leqslant m-1}\gamma_{n+a}\gamma_{m-a-1}(q^{2a}+q^{-2(a+1)})q^{n+m}z^n w^m-\gamma_n\gamma_m z^n w^m
\end{align*}
\noindent Hence the equation to be proved is 
\[
X_{n,m}\gamma_n\gamma_m=\sum_{0\leqslant a\leqslant m}\gamma_{n}\gamma_{m}\eta_a-\sum_{0\leqslant a\leqslant m-1}\gamma_{n+a}\gamma_{m-a-1}(q^{2a}+q^{-2(a+1)})-\gamma_n\gamma_m z^n w^m,
\]
which is equivalent to
\begin{equation}
X_{n,m}=1+(1-t^{-1})(1+t^{2n})\sum_{1\leqslant a\leqslant m}\frac{(t^{1+2n};t^2)_{a-1}}{(t^{2+2n};t^2)_{a-1}}\frac{(t^{2m-2a+2};t^2)_{a}}{(t^{2m-2a+1};t^2)_{a}}\frac{t^a}{1-t^{2(n+a)}}\label{eqn:reduced}
\end{equation}
where we set $t=q^2$. It can be proved by induction with respect to $k$ that the summation over $a=m,m-1,\cdot\cdot\cdot,m-k$ yields  
\[
t^{m-k}\frac{(t^{1+2n};t^2)_{m-k-1}}{(t^{2+2n};t^2)_{m-k-1}}\frac{(t^{2k+2};t^2)_{m-k}}{(t^{2k+1};t^2)_{m-k}}\frac{\sum_{j=0}^{k}t^{2j}}{1-t^{2(n+k)}}.
\]
Setting $k=m-1$ we can see that the right hand side of \eqref{eqn:reduced} is equal to $\frac{t^{2m}-t^{2n}}{1-t^{2(n+m)}}$.

\end{document}